\documentclass[10pt,leqno,twoside]{article}

\usepackage[latin2]{inputenc}
\usepackage{amsmath}
\usepackage{amsfonts}
\usepackage{amssymb}
\usepackage{makeidx}

\newtheorem{theorem}{Theorem}

\newtheorem{lemma}{Lemma}
\newtheorem{example}{Example}
\newtheorem{proposition}{Proposition}
\newtheorem{remark}[theorem]{Remark}

\newfont{\eurorm}{eurm10 scaled 1100}
\newfont{\eurorms}{eurm10 scaled 800}

\def\nfrac#1#2{\mbox{\footnotesize $\displaystyle\frac{#1}{#2}$}}

\usepackage[T1]{fontenc}                         
\usepackage{amsmath,amsfonts}                    
\DeclareMathAlphabet{\mathbbmsl}{U}{bbm}{m}{sl}  

\date{August 2, 2019}

\begin{document}

%
%
%
\title{\bf Some new estimates of precision of Cusa-Huygens and Huygens approximations}
\maketitle

\begin{center}
Branko Male\v sevi\' c${}^{\,\mbox{\scriptsize $1)$}}\!$,
Marija Nenezi\' c${}^{\,\mbox{\scriptsize $1)$}}$,
Ling Zhu${}^{\,\mbox{\scriptsize $\ast,\!2)$}}\!$, \\[0.5 ex]
Bojan Banjac${}^{\,\mbox{\scriptsize $3)$}}$ and Maja Petrovi\' c${}^{\,\mbox{\scriptsize $4)$}}$
\end{center}

\begin{center}
{\footnotesize
\textit{
${}^{\;\mbox{\scriptsize $1)$}}$School of Electrical Engineering, University of Belgrade,                      \\[-0.25 ex]
Bulevar Kralja Aleksandra 73, 11000 Belgrade, Serbia}                                                          \\[0.25 ex]
\textit{
${}^{\;\mbox{\scriptsize $2)$}}$Department of Mathematics, Zhejiang Gongshang University,                      \\[-0.25 ex]
Hangzhou City, Zhejiang Province, 310018, China}                                                               \\[0.25 ex]
\textit{
${}^{\;\mbox{\scriptsize $3)$}}$Faculty of Technical Sciences, University of Novi Sad,                         \\[-0.25 ex]
Trg Dositeja Obradovi\' ca 6, 21000 Novi Sad, Serbia}                                                          \\[0.25 ex]
\textit{
${}^{\;\mbox{\scriptsize $4)$}}$The Faculty of Transport and Traffic Engineering, University of Belgrade,      \\[-0.25 ex]
Vojvode Stepe 305, 11000 Belgrade, Serbia}                                                                     \\[0.25 ex]
}

\end{center}

\medskip
\noindent {\small \textbf{Abstract.} {\small In this paper we present some new upper bounds of the {\sc Cusa}-{\sc Huygens}
and the {\sc Huygens} approximations. Bounds are obtained in the forms of some polynomial~and~some rational functions.}

\footnote{$\!\!\!\!\!\!\!\!\!\!\!\!\!\!$
{\scriptsize
${}^{\mbox{\scriptsize $\ast$}}$Corresponding author.          \\[0.0 ex]
Emails:                                                        \\[0.0 ex]
{\em Branko Male\v sevi\' c} {\tt $<$branko.malesevic@etf.rs$>$},
{\em Marija Nenezi\' c} {\tt $<$maria.nenezic@gmail.com$>$},
{\em Ling Zhu} {\tt $<$zhuling0571@163.com$>$},
{\em Bojan Banjac} {\tt $<$bojan.banjac@uns.ac.rs$>$}, {\em Maja Petrovi\' c} {\tt $<$majapet@sf.bg.ac.rs$>$}}}


{\small \tt MSC: 42A10, 26D05}

\section{Introduction}

In this paper is considered the following {\sc Cusa}-{\sc Huygens} inequality
\begin{equation}
\label{Cusa-Huygens}
\frac{3\sin x}{2+\cos x}
<
x
<
\displaystyle\frac{2}{3} \sin x + \displaystyle\frac{1}{3} \tan x,
\end{equation}
for $x \!\in\! \left( 0, \nfrac{\pi}{2} \right)$, as shown in \cite{D_S_Mitrinovic_1970}, \cite{G_V_Milovanovic_2014}
and \cite{Sandor_Bencze_2005}. Let us emphasize that the following approximation$:$
\begin{equation}
x
\approx
\frac{3\sin x}{2+ \cos x},
\end{equation}
for $x \!\in\! \left(0, \pi \right]$, was first surmised in the {\sc De Cusa}'s {\em Opera book}, see \cite{QUERIES - REPLIES 1949} and \cite{K_T_Vahlen_1911}.
Approximation stated above will be called the {\it Cusa-Huygens approximation}.

\smallskip
Let us consider {\em the error of the {\sc Cusa}-{\sc Huygens} approximation} as the following function$:$
\begin{equation}
R(x) =  x - \frac{3\sin x}{2+\cos x},
\end{equation}
for $x \in [0,\pi]$. One estimation of the precision of the {\sc Cusa}-{\sc Huygens} approximation is given by the following statement of {\sc Ling Zhu}:
\begin{theorem} {\rm \cite{Ling_Zhu_2018}}
It is true that$:$
\begin{equation}
\label{(4)}
\frac{1}{180} x^5 < x - \frac{3\sin x}{2+\cos x},
\end{equation}
and
\begin{equation}
\label{(5)}
\frac{1}{2100} x^7 < x - \frac{3\sin x}{2+\cos x}\!\left(\! 1 + \displaystyle\frac{(1-\cos x)^2}{9(3+2\cos x)} \!\right),
\end{equation}
for $x \!\in\! \left( 0, \pi \right]$ and that $1/180$ and $1/2100$ are the best constants in the previous inequalities, respectively.
\end{theorem}
The results of the previous theorem are corrections of the Theorem 3.4.20 from monograph \cite{D_S_Mitrinovic_1970}.
This important discovery and the resulting corrections took place in 2018, almost half a century after the publication of classics
\cite{Ling_Zhu_2018}.

\smallskip
In this paper we consider also the following \textit{the \textsc{Huygens}'s approximation}$:$
\begin{equation}
\label{Huygens}
x
\approx
\displaystyle\frac{2}{3} \sin x + \displaystyle\frac{1}{3} \tan x,
\end{equation}
for $x \in \left(0, \mbox{\small $\displaystyle\frac{\pi}{2}$} \right)$.
Estimates of \textit{the function of error of {\sc Huygens} approximation}
$
Q(x) = \mbox{\footnotesize $\displaystyle\frac{2}{3}$} \sin x + \mbox{\footnotesize $\displaystyle\frac{1}{3}$}\tan x - x,
$
for $x \!\in\! \left(0, \mbox{\small $\displaystyle\frac{\pi}{2}$} \right)$, are achieved by use of some poly\-nomial functions
and some rational functions. Necessary theoretical basis for that research are stated in the following section.

\section{Preliminaries}

\textbf{Double sided Taylor approximations}\\[0.5 ex]
Let us introduce some denotations and the basic claims that shall be used according to the papers
\cite{B_Malesevic_M_Rasajski_T_Lutovac_2019} and \cite{B_Malesevic_T_Lutovac_M_Rasajski_B_Banjac_2019}.
Let us begin from real function $f : (a, b) \longrightarrow \mathbb{R}$
for which there are the finite values
$f^{(k)}(a+)=\!\lim\limits_{x \rightarrow a+}{f^{(k)}(x)}$, $k=0,1,\ldots,n$, for $n \!\in\! \mathbb{N}_{0}$.
Here we use denotation $T_n^{f,\,a+}(x)$ for {\sc Taylor} polynomial of order $n$, for $n \!\in\! \mathbb{N}_{0}$,
for function $f(x)$ defined in right neighbourhood of $a$:
\begin{equation}
T_n^{f,\,a+}(x)=\displaystyle\sum_{k=0}^{n}{\displaystyle\frac{f^{(k)}(a+)}{k!}(x-a)^k}.
\end{equation}
We shall call $T_n^{f,\,a+}(x)$ {\em the first {\sc Taylor} approximation in the right neighbourhood of~$a$} \cite{B_Malesevic_M_Rasajski_T_Lutovac_2019}.
For \mbox{$n \!\in\! \mathbb{N}_{0}$}, let us define {\em the remainder of first {\sc Taylor} approximation in the right neighbourhood of~$a$} by
$
R_{n}^{f,\,a+}(x)
=
f(x) - T_{n}^{f,\,a+}(x).
$
In the paper \cite{B_Malesevic_M_Rasajski_T_Lutovac_2019} are considered polynomials:
\begin{equation}
\mbox{$\mathbbmsl{T}$}_n^{f;\,a+,\,b-}(x)
=
\left\{
\begin{array}{ccl}
T_{n-1}^{f,\,a+}(x)
+
\displaystyle\frac{1}{(b - a)^n}R_{n-1}^{f,\,a+}(b-)(x-a)^n                     &:& n \geq 1 \\[2.5 ex]
f(b-)                                                                           &:& n = 0,
\end{array}
\right.
\end{equation}
and determined as {\em second {\sc Taylor} approximation in right neighbourhood of $a$}, for $n \!\in\! \mathbb{N}_{0}$,
\cite{B_Malesevic_M_Rasajski_T_Lutovac_2019}. Then following statement is true.

\begin{theorem}
\label{Theorem_1}
Let us assume that $f(x)$ is an real function over $(a,b)$, and that $n$ is natural number such that there exist $f ^{(k)}(a+)$,
for $k \!\in\! \{0,1,2, \ldots ,n\}$. Let us assume that $f^{(n)}(x)$ is increasing over $(a,b)$. Then for every $x \!\in\! (a,b)$
following inequality is true$\,:$
\begin{equation}
\label{(2)}
T_n^{f,\,a+}(x) < f(x) < \mbox{$\mathbbmsl{T}$}_n^{f;\,a+,\,b-}(x).
\end{equation}
At that, if $f^{(n)}(x)$ is decreasing over $(a,b)$, then reversed inequality from \mbox{\rm (\ref{(2)})} is true.
\end{theorem}
Previous statement we call the {\em Theorem on double-sided \textsc{Taylor}'s approximations} in
\cite{B_Malesevic_M_Rasajski_T_Lutovac_2019} and
\cite{B_Malesevic_T_Lutovac_M_Rasajski_B_Banjac_2019}, i.e.$\;${\em Theorem WD} in
\cite{B_Malesevic_T_Lutovac_M_Rasajski_C_Mortici_Adv._Difference_Equ._2018}-\cite{M_Nenenzic_L_Zhu_AADM_2018}.
Let us emphasize that proof of this Theorem (i.e.$\;$Theorem~2 u \cite{S_Wu_L_Debnath_2009}) is based on
{\sc L'Hospital}'s rule for the monotonicity. Similar method is used in proving some close theorems in
\cite{S_Wu_HM_Srivastva_2008a}, \cite{S_Wu_L_Debnath_2008} and \cite{S_Wu_HM_Srivastva_2008b},
which were previously published. Further, the following claims are true.
\begin{proposition} {\rm \cite{B_Malesevic_M_Rasajski_T_Lutovac_2019}}
\label{Proposition_1}
Let $f : (a, b) \longrightarrow \mathbb{R}$ be such real function that there exist the first and the second {\sc Taylor} approximation in the right neighbourhood of $a$, for some $n \in {N}_0$. Then,
\begin{equation}
\label{sgn}
{\mathop{\rm sgn}} {\Big (} \mbox{$\mathbbmsl{T}$}_{n}^{f,\, a+, \, b-}(x)\,-\,\mbox{$\mathbbmsl{T}$}_{n+1}^{f,\, a+, \, b-}(x) {\Big )}
=\,
{\mathop{\rm sgn}} {\Big (} {f(b-)\,-\,T_{n}^{f,\,a+}(b)} {\Big )},
\end{equation}
for every $x \in (a,b)$.
\end{proposition}
\begin{theorem} {\rm \cite{B_Malesevic_M_Rasajski_T_Lutovac_2019}}
\label{Theorem_3}
Let  $f : (a, b) \longrightarrow \mathbb{R}$ be real analytic function with the power series$:$
\begin{equation}
\label{f(x)}
f(x) = \sum_{k=0}^{\infty}{c_k(x-a)^k},
\end{equation}
where $c_k \in \mathbb{R}$ and $c_k \geq 0$ for every $k \in \mathbb{N}_0$. Then,
\begin{equation}
\label{Double_inequalities}
 \begin{array}{c}
 T_0^{f,\,a+}(x) \leq \ldots\leq T_n^{f,\,a+}(x) \leq T_{n+1}^{f,\,a+}(x) \leq \ldots             \\[1.50 ex]
 \ldots \leq  f(x) \leq \ldots                                                                    \\[1.50 ex]
 \ldots \leq  \mbox{$\mathbbmsl{T}$}_{n+1}^{f;\,a+,\,b-}(x)
 \leq   \mbox{$\mathbbmsl{T}$}_n^{f;\,a+,\,b-}(x)
 \leq  \ldots  \leq \mbox{$\mathbbmsl{T}$}_0^{f;\,a+,\,b-}(x),
 \end{array}
\end{equation}
for every $x \in (a,b)$. If $c_k \in \mathbb{R}$ and $c_k \leq 0\,$ for every $k \in \mathbb{N}_0$,
then the reversed inequality is true.
\end{theorem}

\smallskip\noindent
\textbf{Inequality for \textsc{Bernoulli} numbers}\\[0.5 ex]
Let $(\mbox{\bf B}_k)$ be the sequence of {\sc Bernoulli} numbers as it is usually considered, for example see \cite{I_Gradshteyn_I_Ryzhik_2014}.
In this paper we state well the known inequality for {\sc Bernoulli} numbers as given by {\sc D. A'niello}~in~\cite{C_D'Aniello_1994}:
\begin{equation}
\mbox{\small $\displaystyle\frac{2(2n)!}{\pi^{2n}}$}
\mbox{\small $\displaystyle\frac{1}{2^{2n}-1}$}
<
\mbox{\bf B}_{2n}
<
\mbox{\small $\displaystyle\frac{2(2n)!}{\pi^{2n}}$}
\mbox{\small $\displaystyle\frac{1}{2^{2n}-2}$}.
\end{equation}
Previous inequality can be rewritten in the equivalent form
\begin{equation}
\label{Bernoulli's_double_inequality}
2 \mbox{\small $\displaystyle\frac{2^{2n}}{\pi^{2n}}$}
<
\mbox{\small $\displaystyle\frac{2^{2n}(2^{2n}-1)|\mbox{\bf B}_{2n}|}{(2n)!}$}
<
2 \mbox{\small $\displaystyle\frac{2^{2n}}{\pi^{2n}}$} \mbox{\small $\displaystyle\frac{2^{2n}-1}{2^{2n}-2}$},
\end{equation}
for $n \in \mathbb{N}$, and it shall be used in the next section.

\section{Main results}

\medskip
\noindent
{\bf 3.1 Case of Cusa-Huygens approximation}\\[0.5 ex]
In this part we determinate some upper bounds of one estimation of error of the {\sc Cusa}-{\sc Huygens} approximation.

\smallskip
\noindent
In connection with inequality (\ref{(4)}) we consider the following statements.

\begin{lemma}
The function
\begin{equation}
h(t) = \displaystyle\frac{30 \sin t + 15 \cos t \sin t}{4 \cos^2t + 22 \cos t + 19} : [0,\pi] \longrightarrow \mathbb{R}
\end{equation}
has$\,:$

\medskip
\noindent
{\boldmath $1.$}
exactly one maximum on the interval $(0,\pi)$ at the point
\begin{equation}
t_1
=
\pi - \mbox{\rm arccos}
\!\left(\!
1
-
\frac{\sqrt[3]{98+42 \sqrt{105}}}{14}
+
\frac{4}{\sqrt[3]{98+42 \sqrt{105}}}
\right)
=
2.73210 ...
\end{equation}
and the numerical value of the function $h(t)$ in the point $t_1$ is
\begin{equation}
h(t_1)=2.95947 ... \, ;
\end{equation}

\medskip
\noindent
{\boldmath $2.$}
exactly one inflection point on the interval $(0,\pi)$
\begin{equation}
t_2 = \pi - \mbox{\rm arccos} \!\left(\!\displaystyle\frac{35 -
3\sqrt{21}}{28}\!\right) = 2.43258...
\end{equation}
and the numerical value of the function $h(t)$ in the point $t_2$ is
\begin{equation}
h(t_2)=2.63119 ... \, .
\end{equation}

\end{lemma}
\mbox{\bf Proof.} Based on the first and the second derivation of the function $h(t):$
\begin{equation}
h'(t)
=
\displaystyle\frac{210 \cos^3t + 630 \cos^2t + 810\cos t + 375}{(4\cos^2t + 22 \cos t +19)^2}
\end{equation}
and
\begin{equation}
h''(t)
=
\displaystyle\frac{30(28 \cos^2t + 70 \cos t + 37)\,(\cos t - 1)^2 \sin t}{(4 \cos^2t + 22 \cos t + 19)^3},
\end{equation}
statements {\boldmath $1.$} and {\boldmath $2.$} are true. \hfill $\Box$

\begin{lemma}
The equation
\begin{equation}
h(t)=t,
\end{equation}
has exactly one solution
\begin{equation}
t_0 = 2.83982...
\end{equation}
on the interval $(0, \pi)$.
\end{lemma}
{\bf Proof.} The equalities $h(0) = 0$ and $h(\pi) = 0$ are true. The function $h(t)$ is strictly increasing on the interval $(0,t_1)$
and strictly decreasing on the interval $(t_1,\pi)$. The function $h(t)$ is convex on the interval $(0,t_2)$ and concave on the interval $(t_2,\pi)$.
Let us note that
$$
h(t_1) > t_1
\quad \mbox{and} \quad
h(t_2) > t_2.
$$
Therefore we can conclude that in the interval $(t_1,\pi)$ exists exactly one solution of the equation $h(t)=t$ with the numerical value
$t_0=2.83982...\,$.
\hfill $\Box$

\begin{lemma}
The function
\begin{equation}
f(t)
=
\displaystyle\frac{t-\mbox{\footnotesize $\displaystyle\frac{3 \sin t}{2 + \cos t}$}}{t^5} : (0,\pi) \longrightarrow \mathbb{R}
\end{equation}
has exactly one maximum in the point $t_0 = 2.83982...$ and the  numerical value of the function $f(t)$ in the point of the maximum is
\begin{equation}
\mbox{\small $M$}_1 = f(t_0) =  0.010756...\,.
\end{equation}
\end{lemma}
{\bf Proof.} The statement follows from the first derivation
\begin{equation}
f'(t)
=
\displaystyle\frac{30 \sin t + 15 \cos t \sin t - (4 \cos^2 t + 22 \cos t + 19)t}{(2 + \cos t)^2 t^6}
\end{equation}
directly and the previous two lemmas. \hfill $\Box$

\medskip
\noindent
Let us denote
\begin{equation}
m_1 = \frac{1}{\mbox{\small $M$}_1} =  92.96406 ... \, .
\end{equation}
Then, based on the previous three lemmas and result of the paper {\sc Ling Zhu} \cite{Ling_Zhu_2018} we have the following statement.

\begin{theorem}
\label{Theorem_4a}
The following inequalities are true
\begin{equation}
\label{(19)}
\frac{1}{180} x^5 < x - \frac{3\sin x}{2+\cos x} \leq \frac{1}{m_1} x^5 = \frac{1}{92.96406 ...} x^5,
\end{equation}
for $x \!\in\! \left( 0, \pi \right]$.
\end{theorem}
The above consideration on estimates of the precision of the {\sc Cusa}-{\sc Huygens} approximation
may be further generalized if determined  the {\sc Maclaurin} series of the {\sc Cusa}-{\sc Huygens} function$:$
\begin{equation}
\theta(x)=\frac{3 \sin x}{2 + \cos x} : [0,\pi] \longrightarrow \mathbb{R}.
\end{equation}



\break

\noindent
Next, in the connection with the inequality (\ref{(5)}) we consider the following statements.

\begin{lemma}
The function
\begin{equation}
\kappa(\tau) = \displaystyle\frac{294 \sin \tau + 217 \cos \tau \sin
\tau + 14 \cos^2 \tau \sin \tau}{2 \cos^3 \tau + 78 \cos^2 \tau +
258 \cos \tau + 187} : [0,\pi] \longrightarrow \mathbb{R}
\end{equation}
has

\medskip
\noindent
{\boldmath $1.$}
exactly one maximum on the interval $(0,\pi)$ at the point
\begin{equation}
\tau_1
=
2.79340 ...
\end{equation}
and the numerical value of the function $\kappa$ in the point $\tau_1$ is
\begin{equation}
\kappa(\tau_1)=2.97564 ... \, ;
\end{equation}

\medskip
\noindent
{\boldmath $2.$}
exactly one inflection point on the interval $(0,\pi)$
\begin{equation}
\tau_2
=
2.55459...
\end{equation}
and the numerical value of the function $\kappa$ in the point $\tau_2$ is
\begin{equation}
\kappa(\tau_2)=2.71423... \, .
\end{equation}

\end{lemma}
\mbox{\bf Proof.} Based on the first and the second derivation of the function $\kappa(\tau)$
\begin{equation}
\kappa'(\tau) = \nfrac{658 \cos^5 \tau \!+\! 6076 \cos^4 \tau \!+\! 41776 \cos^3 \tau \!+\! 96236 \cos^2 \tau \!+\! 95606 \cos \tau \!+\!35273}{
(2 \cos^3 \tau \!+\! 78 \cos^2 \tau \!+\! 258 \cos \tau \!+\! 187)^2}
\end{equation}
and
\begin{equation}
\!\!
\kappa''(\tau) = \nfrac{14(94 \cos^4 \tau \!-\! 1648 \cos^3 \tau \!-\! 23700 \cos^2 \tau \!-\! 46207 \cos \tau \!-\! 23039)\,
(\cos \tau \!-\! 1)^3 \sin \tau}{(2 \cos^3 \tau \!+\! 78 \cos^2 \tau \!+\! 258 \cos \tau \!+\! 187)^3}
\end{equation}
statements {\boldmath $1.$} and {\boldmath $2.$} are true. \hfill $\Box$

\begin{lemma}
The equation
\begin{equation}
\kappa(\tau)=\tau,
\end{equation}
has exactly one solution
\begin{equation}
\tau_0 = 2.87934...
\end{equation}
on the interval $(0, \pi)$.
\end{lemma}
{\bf Proof.} The equalities $\kappa(0) = 0$ and $\kappa(\pi) = 0$ are true. The function $\kappa(\tau)$ is strictly increasing on the interval
$(0,\tau_1)$ and strictly decreasing on the interval $(\tau_1,\pi)$. The function $\kappa(\tau)$ is convex on the interval $(0,\tau_2)$ and concave
on the interval $(\tau_2,\pi)$. Let us note that
$$
\kappa(\tau_1) > \tau_1 \quad \mbox{and} \quad \kappa(\tau_2) > \tau_2.
$$
Therefore we can conclude that in the interval $(\tau_1,\pi)$ exists exactly one solution of the equation $\kappa(\tau)=\tau$ with the
numerical value $\tau_0=2.87934...\,$. \hfill $\Box$

\begin{lemma}
The function
\begin{equation}
g(\tau) = \displaystyle\frac{\tau-\mbox{\footnotesize
$\displaystyle \frac{3\sin \tau}{2+\cos \tau}\!\left(\! 1 +
\displaystyle\frac{(1-\cos \tau)^2}{9(3+2\cos \tau)}
\!\right)$}}{\tau^7} : (0,\pi) \longrightarrow \mathbb{R}
\end{equation}
has exactly one maximum in the point $\tau_0 = 2.87934...$ and
the numerical value of the function $g(\tau)$ in the point of the maximum is
\begin{equation}
\mbox{\small $M$}_2 = g(\tau_0) =  0.001112...
\end{equation}
\end{lemma}
{\bf Proof.} Statement follows from the first derivation
\begin{equation}
\!\!\!\!\!
g'(\tau) = \nfrac{294 \sin \tau + \cos \tau \sin \tau
(217+ 14\cos \tau) - (2 \cos^3 \tau + 78 \cos^2 \tau + 258 \cos
\tau + 187)\tau}{3(3 + 2\cos \tau)^2 \tau^8}
\end{equation}
directly and the previous two lemmas. \hfill $\Box$

\medskip
\noindent
Let us denote
\begin{equation}
m_2 = \frac{1}{\mbox{\small $M$}_2} =  899.04062... \, .
\end{equation}
Then, based on the previous three lemmas and result of the paper {\sc Ling Zhu} \cite{Ling_Zhu_2018} we have the following statement.

\begin{theorem}
\label{Theorem_4b}
The following inequalities are true
\begin{equation}
\label{(37)}
\frac{1}{2100} x^7 < x - \frac{3\sin x}{2+\cos x}\!\left(\! 1 + \displaystyle\frac{(1-\cos x)^2}{9(3+2\cos x)} \!\right) \leq \frac{1}{m_2} x^7 = \frac{1}{899.04062...} x^7,
\end{equation}
for $x \!\in\! \left( 0, \pi \right]$.
\end{theorem}
The above consideration may be further generalized if the determined {\sc Maclaurin} series of the function$:$
\begin{equation}
\Theta(x) = \frac{3\sin x}{2+\cos x}\!\left(\! 1 + \displaystyle\frac{(1-\cos x)^2}{9(3+2\cos x)} \!\right) : [0,\pi] \longrightarrow \mathbb{R}.
\end{equation}

\bigskip
\noindent
{\bf 3.2 Case of Huygens approximation}\\[0.5 ex]
In this part we determinate some upper bounds of one estimation of the error of the {\sc Huygens} approximation.
Results stated in preliminarily section are applied to the function$:$
\begin{equation}
\varphi(x) = \frac{2}{3} \sin x + \frac{1}{3} \tan x : \left(0, \mbox{\small $\displaystyle\frac{\pi}{2}$}\right) \longrightarrow \mathbb{R},
\end{equation}
for which we shall use the term  \textit{the \textsc{Huygens} function}.

\medskip
\noindent
\textbf{Some polynomial bounds of the \textsc{Huygens} function.}
Let us start from the well known power series from monograph \cite{I_Gradshteyn_I_Ryzhik_2014}
\begin{equation}
\sin x
=
\displaystyle\sum_{k=0}^{\infty}{\mbox{\small $\displaystyle\frac{(-1)^k}{(2k+1)!}$} x^{2k+1}},
\end{equation}
where $x \in \mathbb{R}$ and
\begin{equation}
\tan x
=
\displaystyle\sum_{k=0}^{\infty}{\mbox{\small $\displaystyle\frac{2^{2k+2}(2^{2k+2}-1)|\mbox{\bf B}_{2k+2}|}{(2k+2)!}$} x^{2k+1}},
\end{equation}
where $|x| < \mbox{\small $\displaystyle\frac{\pi}{2}$}$.
Based on the previous two power series follows
\begin{equation}
\varphi(x)
=
x
+
\mbox{\small $\displaystyle\frac{1}{20}$} x^{5}
+
\mbox{\small $\displaystyle\frac{1}{56}$} x^{7}
+
\mbox{\small $\displaystyle\frac{7}{960}$} x^{9}
+
\mbox{\small $\displaystyle\frac{3931}{1330560}$} x^{11}
+
\ldots
=
\displaystyle\sum_{k=0}^{\infty}{a_k x^{2k+1}},
\end{equation}
with coefficients
\begin{equation}
\label{Coeff_a}
a_k
=
\mbox{\small $\displaystyle\frac{2(-1)^k}{3(2k+1)!}$}
+
\mbox{\small $\displaystyle\frac{2^{2k+2}(2^{2k+2}-1)|\mbox{\bf B}_{2k+2}|}{3(2k+2)!}$},
\end{equation}
where $x \in \left( 0, \mbox{\small $\displaystyle\frac{\pi}{2}$} \right)$ and $k \in \mathbb{N}_0$.
Based on the inequalities (\ref{Bernoulli's_double_inequality}) follows that
\begin{equation}
a_k > 0,
\end{equation}
for $k \in \mathbb{N}_0$. From there comes that, based on Theorem \ref{Theorem_3}, the following claim about some polynomial inequalities
for the {\sc Huygens} function is true.
\begin{theorem}
\label{Theorem_4c}
Let there be given the function
$
\varphi(x) = \displaystyle\sum_{k=0}^{\infty}{a_k x^{2k+1}} : \left(0, \mbox{\small $\displaystyle\frac{\pi}{2}$}\right) \longrightarrow \mathbb{R}
$,
with coefficients $(a_k)$ determined with $(\ref{Coeff_a})$. Let it be that $c \!\in\! \left(0, \mbox{\small $\displaystyle\frac{\pi}{2}$}\right)$
is fixed. Then for $x \!\in\! \left(0,c\right)$ holds$:$
\begin{equation}
 \begin{array}{c}
 T_0^{\varphi,\,0+}(x) < \ldots < T_n^{\varphi,\,0+}(x) < T_{n+1}^{\varphi,\,0+}(x) < \ldots               \\[1.0 ex]
 \ldots <  \displaystyle\frac{2}{3} \sin x + \displaystyle\frac{1}{3} \tan x < \ldots       \\[2.0 ex]
 \ldots <  \mbox{$\mathbbmsl{T}$}_{n+1}^{\varphi;\,0+,\,c-}(x)
 <  \mbox{$\mathbbmsl{T}$}_{n}^{\varphi;\,0+,\,c-}(x)
 <  \ldots  < \mbox{$\mathbbmsl{T}$}_0^{\varphi;\,0+,\,c-}(x).
 \end{array}
\end{equation}
\end{theorem}

\begin{example}
Let us introduce some examples of the inequalities obtained for integers $n=0,1,2,3,4,5$.

\medskip
\noindent
{\boldmath $n=0\!:$}
Let it be that $c \!\in\! \left(0, \mbox{\small $\displaystyle\frac{\pi}{2}$}\right)$
is fixed. Then for $x \!\in\! \left(0, c\right)$ holds$:$
\begin{equation}
0 = T_0^{\varphi,\,0+}(x)
<
\frac{2}{3} \sin x + \frac{1}{3} \tan x
<
\mbox{$\mathbbmsl{T}$}_0^{\varphi;\,0+,\,c-}(x) = \frac{2}{3} \sin c + \frac{1}{3} \tan c.
\end{equation}

\medskip
\noindent
{\boldmath $n=1\!:$}
Let it be that $c \!\in\! \left(0, \mbox{\small $\displaystyle\frac{\pi}{2}$}\right)$
is fixed. Then for $x \!\in\! \left(0, c\right)$ holds$:$
\begin{equation}
x = T_1^{\varphi,\,0+}(x)
<
\frac{2}{3} \sin x + \frac{1}{3} \tan x
<
\mbox{$\mathbbmsl{T}$}_1^{\varphi;\,0+,\,c-}(x)
=
\displaystyle\frac{\mbox{\scriptsize $\displaystyle\frac{2}{3}$} \sin c + \mbox{\scriptsize $\displaystyle\frac{1}{3}$} \tan c}{c} x.
\end{equation}

\medskip
\noindent
{\boldmath $n=2,3,4\!:$}
Let it be that $c \!\in\! \left(0, \mbox{\small $\displaystyle\frac{\pi}{2}$}\right)$
is fixed. Then for $x \!\in\! \left(0, c\right)$ holds$:$
\begin{equation}
x = T_n^{\varphi,\,0+}(x)
<
\frac{2}{3} \sin x + \frac{1}{3} \tan x
<
\mbox{$\mathbbmsl{T}$}_n^{\varphi;\,0+,\,c-}(x)
=
x + \displaystyle\frac{\mbox{\scriptsize $\displaystyle\frac{2}{3}$} \sin c + \mbox{\scriptsize $\displaystyle\frac{1}{3}$} \tan c}{c^n} x^n,
\end{equation}
for integers $n=2,3,4$.

\break

\noindent
{\boldmath $n=5\!:$}
Let it be that $c \!\in\! \left(0, \mbox{\small $\displaystyle\frac{\pi}{2}$}\right)$
is fixed. Then for $x \!\in\! \left(0, c\right)$ holds$:$
\begin{equation}
x + \mbox{\small $\displaystyle\frac{1}{20}$}x^5 = T_5^{\varphi,\,0+}(x)
<
\frac{2}{3} \sin x + \frac{1}{3} \tan x
<
\mbox{$\mathbbmsl{T}$}_5^{\varphi;\,0+,\,c-}(x)
=
\displaystyle\frac{\mbox{\scriptsize $\displaystyle\frac{2}{3}$} \sin c + \mbox{\scriptsize $\displaystyle\frac{1}{3}$} \tan c}{c^5} x^5.
\end{equation}
\end{example}

\smallskip\noindent
From the previous Theorem directly follows
estimate of the function of error of {\sc Huygens} approximation
$
Q(x) = \mbox{\footnotesize $\displaystyle\frac{2}{3}$} \sin x + \mbox{\footnotesize $\displaystyle\frac{1}{3}$} \tan x - x
$
with previously considered polynomial functions.
\begin{theorem}
\label{Theorem_4d}
Let it be given the function
$
\varphi(x) = \displaystyle\sum_{k=0}^{\infty}{a_k x^{2k+1}} : \left(0, \mbox{\small $\displaystyle\frac{\pi}{2}$}\right) \longrightarrow \mathbb{R}
$,
with coefficients  $(a_k)$ determined with $(\ref{Coeff_a})$. Let it be that $c \!\in\! \left(0, \mbox{\small $\displaystyle\frac{\pi}{2}$}\right)$
is fixed. Then for $x \!\in\! \left(0, c\right)$ holds$:$
\begin{equation}
 \begin{array}{c}
 T_0^{\varphi,\,0+}(x) - x < \ldots < T_n^{\varphi,\,0+}(x) - x < T_{n+1}^{\varphi,\,0+}(x) - x < \ldots            \\[1.0 ex]
 \ldots <  \displaystyle\frac{2}{3} \sin x + \displaystyle\frac{1}{3} \tan x - x < \ldots         \\[2.0 ex]
 \ldots <  \mbox{$\mathbbmsl{T}$}_{n+1}^{\varphi;\,0+,\,c-}(x) - x
 <  \mbox{$\mathbbmsl{T}$}_{n}^{\varphi;\,0+,\,c-}(x) - x
 <  \ldots  < \mbox{$\mathbbmsl{T}$}_0^{\varphi;\,0+,\,c-}(x) - x.
 \end{array}
\end{equation}
\end{theorem}

\bigskip
\noindent
\textbf{Some rational bounds for the \textsc{Huygens} function.}
In this section are considered some series for tangent function obtained based on the
well known series for cotangent function \cite{I_Gradshteyn_I_Ryzhik_2014}$:$
\begin{equation}
\begin{array}{rcl}
\mbox{\rm cot}\,x
\!&\!\!=\!\!&\!
\mbox{\small $\displaystyle\frac{1}{x}$}
-
\displaystyle\sum_{k=0}^{\infty}{\mbox{\small $\displaystyle\frac{2^{2k+2}|\mbox{\bf B}_{2k+2}|}{(2k+2)!}$}x^{2k+1}}    \\[3.0 ex]
\!&\!\!=\!\!&\!
\mbox{\small $\displaystyle\frac{1}{x}$}
-
\mbox{\small $\displaystyle\frac{1}{3}$}x
-
\mbox{\small $\displaystyle\frac{1}{45}$}x^3
-
\mbox{\small $\displaystyle\frac{2}{945}$}x^5
-
\ldots
\end{array}
\end{equation}
which converges for $0 < |x| < \pi$. From previous series we conclude that
\begin{equation}
\begin{array}{rcl}
\mbox{\rm tan}\,x
\!&\!\!=\!\!&\!
\mbox{\small $\displaystyle\frac{1}{\mbox{\footnotesize $\displaystyle\frac{\pi}{2}$} - x }$}
-
\displaystyle\sum_{k=0}^{\infty}{\mbox{\small $\displaystyle\frac{2^{2k+2}|\mbox{\bf B}_{2k+2}|}{(2k+2)!}$}
\!\left(\mbox{\small $\displaystyle\frac{\pi}{2}$} - x \right)^{2k+1}}                                                  \\[4.0 ex]
\!&\!\!=\!\!&\!
\mbox{\small $\displaystyle\frac{1}{\mbox{\footnotesize $\displaystyle\frac{\pi}{2}$} - x}$}
-
\mbox{\small $\displaystyle\frac{1}{3}$}\!\left(\mbox{\small $\displaystyle\frac{\pi}{2}$} - x \right)
-
\mbox{\small $\displaystyle\frac{1}{45}$}\!\left(\mbox{\small $\displaystyle\frac{\pi}{2}$} - x \right)^{\!3}
-
\mbox{\small $\displaystyle\frac{2}{945}$}\!\left(\mbox{\small $\displaystyle\frac{\pi}{2}$} - x \right)^{\!5}
-
\ldots
\end{array}
\end{equation}
for $0 \!<\!\left|\mbox{\small $\displaystyle\frac{\pi}{2}$} \!-\! x \right| \!<\! \pi$,
which holds for  $x \!\in\! \left(0, \mbox{\small $\displaystyle\frac{\pi}{2}$}\right)$.
From there
\mbox{$
\phi(x) = \tan x - \mbox{\small $\displaystyle\frac{1}{\mbox{\footnotesize $\displaystyle\frac{\pi}{2}$} - x}$}
$}
determines real analytic function on $\left(0,\mbox{\small $\displaystyle\frac{\pi}{2}$}\right)$.
Let us notice that  {\sc M. Nenezi\' c} and {\sc L. Zhu} in the paper \cite{M_Nenenzic_L_Zhu_AADM_2018}
obtained the following series$:$
\begin{equation}
\label{Razvoj_2}
\!\!\!\!\!\!
\begin{array}{rcl}
\tan x
\!\!&\!\!=\!\!&\!\!
\mbox{\small $\displaystyle\frac{1}{\mbox{\footnotesize $\displaystyle\frac{\pi}{2}$} - x }$}
+
\displaystyle\sum_{k=0}^{\infty}{b_k x^{k}},                                                      \\[4.0 ex]
\!\!&\!\!=\!\!&\!\!
\mbox{\small $\displaystyle\frac{1}{\mbox{\footnotesize $\displaystyle\frac{\pi}{2}$} - x }$}
-
\mbox{\footnotesize $\displaystyle\frac{2}{\pi}$}
+
\left( 1 \!-\! \mbox{\footnotesize $\displaystyle\frac{4}{\pi^2}$} \right)\!x
-
\mbox{\footnotesize $\displaystyle\frac{8}{\pi^3}$} x^2
+
\left( \mbox{\footnotesize $\displaystyle\frac{1}{3}$} \!-\! \mbox{\footnotesize $\displaystyle\frac{16}{\pi^4}$} \right)\!x^3
-
\mbox{\footnotesize $\displaystyle\frac{32}{\pi^5}$} x^4
+
\ldots
\end{array}
\end{equation}
for $x \in \left(0, \mbox{\small $\displaystyle\frac{\pi}{2}$} \right)$, with coefficients
\begin{equation}
b_{k}
=
\left\{
\begin{array}{ccc}
-
\mbox{\small $\displaystyle\frac{2}{\pi}$} &:& k=0                      \\[2.0 ex]
\mbox{\small $\displaystyle\frac{2^{k+1}(2^{k+1}-1)|\mbox{\bf B}_{k+1}|}{(k+1)!}$}
-
\mbox{\small $\displaystyle\frac{2^{k+1}}{\pi^{k+1}}$} &:& k>0,
\end{array}
\right.
\end{equation}
for $k \!\in\! \mathbb{N}_{0}$. Let us introduce sequence
\begin{equation}
\label{Coeff_beta}
\beta_k
=
(-1)^{k-1} b_k
=
\left\{
\begin{array}{ccc}
\mbox{\small $\displaystyle\frac{2^{k+1}}{\pi^{k+1}}$} &:& k = 2 \ell                      \\[2.0 ex]
\mbox{\small $\displaystyle\frac{2^{k+1}(2^{k+1}-1)|\mbox{\bf B}_{k+1}|}{(k+1)!}$}
-
\mbox{\small $\displaystyle\frac{2^{k+1}}{\pi^{k+1}}$} &:& k = 2 \ell - 1,
\end{array}
\right.
\end{equation}
for $\ell \in \mathbb{N}_0$.
Based on the inequality (\ref{Bernoulli's_double_inequality}) the following statement is easily checked.
\begin{lemma}
For fixed $x \!\in\! \left(0,\mbox{\small $\displaystyle\frac{\pi}{2}$} \right)$ and sequence $(\beta_k)$ holds$:$
\begin{equation}
\beta_k x^k > 0, \quad \mathop{\lim}\limits_{k \rightarrow \infty}{\beta_k x^k} = 0, \quad (\beta_k x^k)\downarrow \, .
\end{equation}
\end{lemma}

\medskip
\noindent
Based on the {\sc Leibnitz} alternating series test follows next statement about some rational inequalities for tangent function.
\begin{theorem}
\label{Theorem_5}
Let it be given the function$:$
\begin{equation}
\phi(x)
=
\tan x - \mbox{\small $\displaystyle\frac{1}{\mbox{\footnotesize $\displaystyle\frac{\pi}{2}$} - x}$}
=
\displaystyle\sum_{k=0}^{\infty}{(-1)^{k-1} \beta_k x^{k}} : \left(0, \mbox{\small $\displaystyle\frac{\pi}{2}$}\right)
\longrightarrow \mathbb{R},
\end{equation}
with coefficients $(\beta_k)$ determined by $(\ref{Coeff_beta})$. Then for
$x \in \left(0,\mbox{\small $\displaystyle\frac{\pi}{2}$} \right)$ holds$:$
\begin{equation}
\!\!\!\!\!\!\!\!
\begin{array}{c}
\mbox{\small $\displaystyle\frac{1}{\mbox{\footnotesize $\displaystyle\frac{\pi}{2}$} \!-\! x }$}
+
T_{0}^{\phi,0+}(x)
<
\ldots
<
\mbox{\small $\displaystyle\frac{1}{\mbox{\footnotesize $\displaystyle\frac{\pi}{2}$} \!-\! x }$}
+
T_{2n}^{\phi,0+}(x)
<
\mbox{\small $\displaystyle\frac{1}{\mbox{\footnotesize $\displaystyle\frac{\pi}{2}$} \!-\! x }$}
+
T_{2n+2}^{\phi,0+}(x)
<
\ldots                                                                                            \\[3.5 ex]
<
\tan x
<                                                                                                 \\[1.0 ex]
\ldots
<
\mbox{\small $\displaystyle\frac{1}{\mbox{\footnotesize $\displaystyle\frac{\pi}{2}$} \!-\! x }$}
+
T_{2n+1}^{\phi,0+}(x)
<
\mbox{\small $\displaystyle\frac{1}{\mbox{\footnotesize $\displaystyle\frac{\pi}{2}$} \!-\! x }$}
+
T_{2n-1}^{\phi,0+}(x)
<
\ldots
<
\mbox{\small $\displaystyle\frac{1}{\mbox{\footnotesize $\displaystyle\frac{\pi}{2}$} \!-\! x }$}
+
T_{1}^{\phi,0+}(x).
\end{array}
\!\!\!\!\!\!\!\!
\end{equation}
\end{theorem}

\medskip
\noindent
Furthermore, let us consider the function
\begin{equation}
\psi(x)
\!=\!
\varphi(x)
\!-\!
\mbox{\small $\displaystyle\frac{1}{3}$}
\mbox{\small $\displaystyle\frac{1}{\mbox{\footnotesize $\displaystyle\frac{\pi}{2}$} \!-\! x }$}
\!=\!
\frac{2}{3} \sin x + \frac{1}{3} \tan x
\!-\!
\mbox{\small $\displaystyle\frac{1}{3}$}
\mbox{\small $\displaystyle\frac{1}{\mbox{\footnotesize $\displaystyle\frac{\pi}{2}$} \!-\! x }$}
\!=\!
\displaystyle\sum_{k=0}^{\infty}{c_k x^{k}} \!:\! \left(0, \mbox{\small $\displaystyle\frac{\pi}{2}$}\right)
\longrightarrow \mathbb{R},
\end{equation}
with coefficients $(c_k)$ determined with
\begin{equation}
\label{Coeff_gamma}
c_k
=
\left\{
\begin{array}{ccl}
-\mbox{\small $\displaystyle\frac{2}{3 \pi}$}                                   &:&               k=0    \\[2.50 ex]
\mbox{\small $\displaystyle\frac{2(-1)^{\frac{k-1}{2}}}{3k!}$}
+
\mbox{\small $\displaystyle\frac{2^{k+1}(2^{k+1}\!-\!1)|\mbox{\bf B}_{k+1}|}{3(k+1)!}$}
-
\mbox{\small $\displaystyle\frac{2^{k+1}}{3 \pi^{k+1}}$}                        &:&               k=2j+1 \\[3.00 ex]
-\mbox{\small $\displaystyle\frac{2^{k+1}}{3 \pi^{k+1}}$}                       &:&               k=2j+2,
\end{array}
\right.
\end{equation}
for $j \!\in\! \mathbb{N}_{0}$. Let us introduce the sequence
\begin{equation}
\gamma_k
\!=\!
(-1)^{k-1} c_k
\!=\!
\left\{\!
\begin{array}{ccl}
\mbox{\small $\displaystyle\frac{2}{3 \pi}$}                                   \!\!&\!:\!&\!\!               k\!=\!0        \\[2.50 ex]
\!\!
\mbox{\small $\displaystyle\frac{2(-1)^{\frac{k-1}{2}}}{3k!}$}
+
\mbox{\small $\displaystyle\frac{2^{k+1}(2^{k+1}\!-\!1)|\mbox{\bf B}_{k+1}|}{3(k+1)!}$}
-
\mbox{\small $\displaystyle\frac{2^{k+1}}{3 \pi^{k+1}}$}
                                                                               \!\!&\!:\!&\!\!               k\!=\!2j\!+\!1 \\[3.0 ex]
\mbox{\small $\displaystyle\frac{2^{k+1}}{3 \pi^{k+1}}$}                       \!\!&\!:\!&\!\!               k\!=\!2j\!+\!2,
\end{array}
\right.
\!\!\!\!
\end{equation}
for $j \!\in\! \mathbb{N}_{0}$. By application of symbolic algebra system we can determine
initial part of the power series of $\psi(x)$, for example up to sixth degree$:$
\begin{equation}
\psi(x)
\!=\!
-\mbox{\small $\displaystyle\frac{2}{3 \pi}$}
\!+\!
\left(\! 1 - \mbox{\small $\displaystyle\frac{4}{3 \pi^2}$} \!\right) \! x
\!-\!
\mbox{\small $\displaystyle\frac{8}{3 \pi^3}$} x^2
\!-\!
\mbox{\small $\displaystyle\frac{16}{3 \pi^4}$} x^3
\!-\!
\mbox{\small $\displaystyle\frac{32}{3 \pi^5}$} x^4
\!+\!
\left(\! \mbox{\small $\displaystyle\frac{1}{20}$} \!-\! \mbox{\small $\displaystyle\frac{64}{3 \pi^6}$} \!\right) \! x^5
\!-\!
\mbox{\small $\displaystyle\frac{128}{3 \pi^7}$} x^6
\!+\!
\ldots
\end{equation}
Based on the inequality (\ref{Bernoulli's_double_inequality}) the following statement is simply checked.
\begin{lemma}
For fixed $x \!\in\! \left(0,\mbox{\small $\displaystyle\frac{\pi}{2}$} \right)$ and sequence $(\gamma_k)$ holds$:$
\begin{equation}
\gamma_k x^k > 0\;(\mbox{for $k > 3$}),
\quad \mathop{\lim}\limits_{k \rightarrow \infty}{\gamma_k x^k} = 0,
\quad (\gamma_k x^k)\downarrow \, .
\end{equation}
\end{lemma}
Based on the {\sc Leibnitz} alternating series test follows statement about some rational inequalities for the {\sc Huygens} function.
\begin{theorem}
\label{Theorem_7}
Let there be given the function$:$
\begin{equation}
\psi(x)
=
\varphi(x)
-
\mbox{\small $\displaystyle\frac{1}{3}$}
\mbox{\small $\displaystyle\frac{1}{\mbox{\footnotesize $\displaystyle\frac{\pi}{2}$} \!-\! x }$}
=
\displaystyle\sum_{k=0}^{\infty}{(-1)^{k-1} \gamma_k x^{k}} : \left(0, \mbox{\small $\displaystyle\frac{\pi}{2}$}\right)
\longrightarrow \mathbb{R},
\end{equation}
with coefficients $(\gamma_k)$ determined by  $(\ref{Coeff_gamma})$.
Then for $x \!\in\! \left(0, \mbox{\small $\displaystyle\frac{\pi}{2}$}\right)$ holds$:$
\begin{equation}
\!\!\!\!\!\!\!
\!\!\!\!\!\!\!
\!\!\!\!\!\!\!
\!\!\!\!\!\!\!
\begin{array}{c}
\mbox{\small $\displaystyle\frac{1}{3}$}
\mbox{\small $\displaystyle\frac{1}{\mbox{\footnotesize $\displaystyle\frac{\pi}{2}$} \!-\! x }$}
+
T_{4}^{\psi,0+}(x)
<
\ldots
<
\mbox{\small $\displaystyle\frac{1}{3}$}
\mbox{\small $\displaystyle\frac{1}{\mbox{\footnotesize $\displaystyle\frac{\pi}{2}$} \!-\! x }$}
+
T_{2n}^{\psi,0+}(x)
<
\mbox{\small $\displaystyle\frac{1}{3}$}
\mbox{\small $\displaystyle\frac{1}{\mbox{\footnotesize $\displaystyle\frac{\pi}{2}$} \!-\! x }$}
+
T_{2n+2}^{\psi,0+}(x)
<
\ldots                                                                                            \\[4.5 ex]
<
\mbox{\small $\displaystyle\frac{2}{3}$} \sin x
+
\mbox{\small $\displaystyle\frac{1}{3}$} \tan x
<                                                                                                 \\[3.0 ex]
\ldots
<
\mbox{\small $\displaystyle\frac{1}{3}$}
\mbox{\small $\displaystyle\frac{1}{\mbox{\footnotesize $\displaystyle\frac{\pi}{2}$} \!-\! x }$}
+
T_{2n+1}^{\psi,0+}(x)
<
\mbox{\small $\displaystyle\frac{1}{3}$}
\mbox{\small $\displaystyle\frac{1}{\mbox{\footnotesize $\displaystyle\frac{\pi}{2}$} \!-\! x }$}
+
T_{2n-1}^{\psi,0+}(x)
<
\ldots
<
\mbox{\small $\displaystyle\frac{1}{3}$}
\mbox{\small $\displaystyle\frac{1}{\mbox{\footnotesize $\displaystyle\frac{\pi}{2}$} \!-\! x }$}
+
T_{3}^{\psi,0+}(x).
\end{array}
\!\!\!\!\!\!\!
\!\!\!\!\!\!\!
\end{equation}
\end{theorem}

\begin{example}
Let us introduce some examples of inequalities obtained for integers $n=2, 3$.

\smallskip\noindent
{\boldmath $n=2\!:$} 
For $x \!\in\! \left(0, \mbox{\small $\displaystyle\frac{\pi}{2}$}\right)$ holds$:$
\begin{equation}
\begin{array}{c}
\mbox{\small $\displaystyle\frac{1}{3}$}
\mbox{\small $\displaystyle\frac{1}{\mbox{\footnotesize $\displaystyle\frac{\pi}{2}$} \!-\! x }$}
-
\mbox{\small $\displaystyle\frac{2}{3 \pi}$}
+
\left(\! 1 - \mbox{\small $\displaystyle\frac{4}{3 \pi^2}$} \!\right) \! x
\!-\!
\mbox{\small $\displaystyle\frac{8}{3 \pi^3}$} x^2
\!-\!
\mbox{\small $\displaystyle\frac{16}{3 \pi^4}$} x^3
\!-\!
\mbox{\small $\displaystyle\frac{32}{3 \pi^5}$} x^4
=                                                                                                 \\[3.0 ex]
=
\mbox{\small $\displaystyle\frac{1}{3}$}
\mbox{\small $\displaystyle\frac{1}{\mbox{\footnotesize $\displaystyle\frac{\pi}{2}$} \!-\! x }$}
+
T_{4}^{\psi,0+}(x)
<                                                                                                 \\[4.0 ex]
<
\displaystyle\frac{2}{3} \sin x + \displaystyle\frac{1}{3} \tan x
<                                                                                                 \\[3.0 ex]
<
\mbox{\small $\displaystyle\frac{1}{3}$}
\mbox{\small $\displaystyle\frac{1}{\mbox{\footnotesize $\displaystyle\frac{\pi}{2}$} \!-\! x }$}
+
T_{3}^{\psi,0+}(x)
=                                                                                                 \\[3.0 ex]
=
\mbox{\small $\displaystyle\frac{1}{3}$}
\mbox{\small $\displaystyle\frac{1}{\mbox{\footnotesize $\displaystyle\frac{\pi}{2}$} \!-\! x }$}
-
\mbox{\small $\displaystyle\frac{2}{3 \pi}$}
+
\left(\! 1 - \mbox{\small $\displaystyle\frac{4}{3 \pi^2}$} \!\right) \! x
\!-\!
\mbox{\small $\displaystyle\frac{8}{3 \pi^3}$} x^2
\!-\!
\mbox{\small $\displaystyle\frac{16}{3 \pi^4}$} x^3.
\end{array}
\end{equation}

\medskip
\noindent
{\boldmath $n=3\!:$} 
For $x \!\in\! \left(0, \mbox{\small $\displaystyle\frac{\pi}{2}$}\right)$ holds$:$
\begin{equation}
\!\!\!\!\!\!\!\!\!\!\!\!\!\!\!\!\!\!\!\!\!
\begin{array}{c}
\mbox{\small $\displaystyle\frac{1}{3}$}
\mbox{\small $\displaystyle\frac{1}{\mbox{\footnotesize $\displaystyle\frac{\pi}{2}$} \!-\! x }$}
-
\mbox{\small $\displaystyle\frac{2}{3 \pi}$}
\!+\!
\left(\! 1 - \mbox{\small $\displaystyle\frac{4}{3 \pi^2}$} \!\right) \! x
\!-\!
\mbox{\small $\displaystyle\frac{8}{3 \pi^3}$} x^2
\!-\!
\mbox{\small $\displaystyle\frac{16}{3 \pi^4}$} x^3
\!-\!
\mbox{\small $\displaystyle\frac{32}{3 \pi^5}$} x^4
\!+\!
\left(\! \mbox{\small $\displaystyle\frac{1}{20}$} \!-\! \mbox{\small $\displaystyle\frac{64}{3 \pi^6}$} \!\right) \! x^5
\!-\!
\mbox{\small $\displaystyle\frac{128}{3 \pi^7}$} x^6                                              \\[3.0 ex]
=
\mbox{\small $\displaystyle\frac{1}{3}$}
\mbox{\small $\displaystyle\frac{1}{\mbox{\footnotesize $\displaystyle\frac{\pi}{2}$} \!-\! x }$}
+
T_{6}^{\psi,0+}(x)
<                                                                                                 \\[4.0 ex]
<
\displaystyle\frac{2}{3} \sin x + \displaystyle\frac{1}{3} \tan x
<                                                                                                 \\[3.0 ex]
<
\mbox{\small $\displaystyle\frac{1}{3}$}
\mbox{\small $\displaystyle\frac{1}{\mbox{\footnotesize $\displaystyle\frac{\pi}{2}$} \!-\! x }$}
+
T_{5}^{\psi,0+}(x)
=                                                                                                 \\[3.0 ex]
=
\mbox{\small $\displaystyle\frac{1}{3}$}
\mbox{\small $\displaystyle\frac{1}{\mbox{\footnotesize $\displaystyle\frac{\pi}{2}$} \!-\! x }$}
-
\mbox{\small $\displaystyle\frac{2}{3 \pi}$}
\!+\!
\left(\! 1 - \mbox{\small $\displaystyle\frac{4}{3 \pi^2}$} \!\right) \! x
\!-\!
\mbox{\small $\displaystyle\frac{8}{3 \pi^3}$} x^2
\!-\!
\mbox{\small $\displaystyle\frac{16}{3 \pi^4}$} x^3
\!-\!
\mbox{\small $\displaystyle\frac{32}{3 \pi^5}$} x^4
\!+\!
\left(\! \mbox{\small $\displaystyle\frac{1}{20}$} \!-\! \mbox{\small $\displaystyle\frac{64}{3 \pi^6}$} \!\right) \! x^5.
\end{array}
\!\!\!\!\!\!\!\!\!\!\!
\end{equation}
\end{example}
\begin{remark}
For $x \!\in\! \left(0, \mbox{\small $\displaystyle\frac{\pi}{2}$}\right)$ holds$:$
\begin{equation}
\displaystyle\frac{2}{3} \sin x + \displaystyle\frac{1}{3} \tan x
<
\mbox{\small $\displaystyle\frac{1}{3}$}
\mbox{\small $\displaystyle\frac{1}{\mbox{\footnotesize $\displaystyle\frac{\pi}{2}$} \!-\! x }$}
+
T_{2}^{\psi,0+}(x)
<
\mbox{\small $\displaystyle\frac{1}{3}$}
\mbox{\small $\displaystyle\frac{1}{\mbox{\footnotesize $\displaystyle\frac{\pi}{2}$} \!-\! x }$}
+
T_{1}^{\psi,0+}(x).
\end{equation}
\end{remark}
From the previous Theorem simply follows
estimate of the function of the error of {\sc Huygens} approximation
$
Q(x) = \mbox{\footnotesize $\displaystyle\frac{2}{3}$} \sin x + \mbox{\footnotesize $\displaystyle\frac{1}{3}$} \tan x - x
$
with previously considered rational functions.
\begin{theorem}
\label{Theorem_7b}
Let it be given the function$:$
\begin{equation}
\psi(x)
=
\phi(x)
-
\mbox{\small $\displaystyle\frac{1}{3}$}
\mbox{\small $\displaystyle\frac{1}{\mbox{\footnotesize $\displaystyle\frac{\pi}{2}$} \!-\! x }$}
=
\displaystyle\sum_{k=0}^{\infty}{(-1)^{k-1} \gamma_k x^{k}} : \left(0, \mbox{\small $\displaystyle\frac{\pi}{2}$}\right)
\longrightarrow \mathbb{R},
\end{equation}
with coefficients $(\gamma_k)$ determined by $(\ref{Coeff_gamma})$.
Then for $x \!\in\! \left(0, \mbox{\small $\displaystyle\frac{\pi}{2}$}\right)$ holds
\begin{equation}
\!\!\!\!\!\!\!\!
\!\!\!\!\!\!\!\!
\!\!\!\!\!\!\!\!
\!\!\!\!\!\!\!\!
\!\!\!\!
\begin{array}{c}
\mbox{\small $\displaystyle\frac{1}{3}$}
\mbox{\small $\displaystyle\frac{1}{\mbox{\footnotesize $\displaystyle\frac{\pi}{2}$} \!-\! x }$}
+
T_{4}^{\psi,0+}(x) - x
<
\ldots
<
\mbox{\small $\displaystyle\frac{1}{3}$}
\mbox{\small $\displaystyle\frac{1}{\mbox{\footnotesize $\displaystyle\frac{\pi}{2}$} \!-\! x }$}
+
T_{2n}^{\psi,0+}(x) - x
<
\mbox{\small $\displaystyle\frac{1}{3}$}
\mbox{\small $\displaystyle\frac{1}{\mbox{\footnotesize $\displaystyle\frac{\pi}{2}$} \!-\! x }$}
+
T_{2n+2}^{\psi,0+}(x) - x
<
\ldots                                                                                            \\[4.5 ex]
<
\mbox{\small $\displaystyle\frac{2}{3}$} \sin x
+
\mbox{\small $\displaystyle\frac{1}{3}$} \tan x - x
<                                                                                                 \\[3.0 ex]
\ldots
<
\mbox{\small $\displaystyle\frac{1}{3}$}
\mbox{\small $\displaystyle\frac{1}{\mbox{\footnotesize $\displaystyle\frac{\pi}{2}$} \!-\! x }$}
+
T_{2n+1}^{\psi,0+}(x) - x
<
\mbox{\small $\displaystyle\frac{1}{3}$}
\mbox{\small $\displaystyle\frac{1}{\mbox{\footnotesize $\displaystyle\frac{\pi}{2}$} \!-\! x }$}
+
T_{2n-1}^{\psi,0+}(x) - x
<
\ldots
<
\mbox{\small $\displaystyle\frac{1}{3}$}
\mbox{\small $\displaystyle\frac{1}{\mbox{\footnotesize $\displaystyle\frac{\pi}{2}$} \!-\! x }$}
+
T_{3}^{\psi,0+}(x) - x.
\end{array}
\!\!\!\!\!\!\!\!
\!\!\!\!\!\!\!\!
\!\!\!\!\!\!\!\!
\!\!\!\!\!\!\!\!
\end{equation}
\end{theorem}

\section{Conclusion}

Based on the inequalities (\ref{(4)}) and (\ref{(5)}), stated by {\sc Zhu} \cite{Ling_Zhu_2018}, using elementary analysis
we have obtained in Theorems \ref{Theorem_4a} and \ref{Theorem_4b}
two new double inequalities which can be used to estimate
some polynomial bounds of the function of error of the {\sc Cusa}-{\sc Huygens} approximation.
With Theorem~\ref{Theorem_4d} were determined some bounds of the function of error of the  {\sc Huygens}
approximation using polynomial functions and with Theorem \ref{Theorem_7b} were determined some bounds of the function of error of the \textsc{Huygens}
approximation using rational functions. Let us emphasize that with the Theorem \ref{Theorem_5} were given some bounds of tangent function
by use of rational functions which can be applied to other parts of Theory of analytical inequalities. Lastly, let us notice that the proofs
of the considered inequalities~can~be~also~obtained by application of methods and algorithms presented in papers \cite{Ling_Zhu_2009},
\cite{F_Qi_DW_Niu_BN_Guo_2009}, \cite{C_Mortici_2011},
\mbox{\cite{B_Malesevic_2007}-\cite{B_Malesevic_T_Lutovac_M_Rasajski_C_Mortici_Adv._Difference_Equ._2018}},
\cite{Chen_Shi_Wang_Xiang_2017}-\cite{Chen_Ma_Li_2019}
and in dissertation \cite{B_D_Banjac_2019}.

\bigskip

\bigskip
\noindent \textbf{Acknowledgement.}
The authors are grateful to the reviewers for their careful reading and for their valuable comments.

\bigskip
\noindent \textbf{Funding.} The research of the first author was supported in part by the Serbian Ministry of Education,
Science and Technological Development, under Projects ON 174032 \& III 44006. The research of the third author was
supported in part by the Natural Science Foundation of China grants No. 61772025.

\bigskip
\noindent \textbf{Competing Interests.}
The authors would like to state that they do not have any competing interests in the subject of this research.

\bigskip
\noindent \textbf{Author's Contributions.} All the authors participated in every phase of the research conducted for this paper.

\end{document}